\documentclass{article}[12pt]

\usepackage{stmaryrd}
\usepackage{times,fancyhdr}
\usepackage{amsmath,
                    graphicx,amssymb}
\usepackage{graphicx}
\usepackage{hyperref}
\usepackage{tikz}
\usetikzlibrary{matrix, positioning}
\input xy
\xyoption{all}
\usepackage{xcolor}
\usepackage{amsthm}
\def\tensors{\mathfrak T^{3} ({\cal L})}
\def\1tensors{\mathfrak T^{3}_\omega ({\cal L})}
\def\2tensors{\mathfrak T^{3}_{\overline\omega}({\cal L})}
\newtheorem{theorem}{Theorem}

\newtheorem{definition}{Definition}

\title{Ternary Associativity and Ternary Lie Algebra at Cube Root of Unity}

\author{Viktor Abramov\\ Institute of Mathematics and Statistics\\ University of Tartu, Estonia}
\date{}

\begin{document}
\maketitle

\abstract{We propose a new approach to extending the notion of commutator and Lie algebra to algebras with ternary multiplication laws. Our approach is based on ternary associativity of the first and second kind. We propose a ternary commutator, which is a linear combination of six (all permutations of three elements) triple products. The coefficients of this linear combination are the cube roots of unity. We find an identity for the ternary commutator that holds due to ternary associativity of the first or second kind. The form of the found identity is determined by the permutations of the general affine group $GA(1,5)\subset S_5$. We consider the found identity as an analogue of the Jacobi identity in the ternary case. We introduce the concept of a ternary Lie algebra at the cubic root of unity and give examples of such an algebra constructed using ternary multiplications of rectangular and three-dimensional matrices. We point out the connection between the structure constants of a ternary Lie algebra with three generators and an irreducible representation of the rotation group.}

\section{Introduction}
The concept of a group endowed with the structure of a smooth manifold and of its tangent space at the identity of a group endowed with the structure of a Lie algebra plays an important role in differential geometry, classical mechanics and theoretical physics. The development of the theory of Lie groups and algebras began with the work of Sophus Lie on the study of symmetries of differential equations and can be considered as an analogue of Galois theory for differential equations. The development of the theory of Lie groups and algebras is closely intertwined with the development of theoretical physics. The development of supersymmetric field theories that emerged in the 1970s is based on the concept of Lie superalgebra, which can be considered as a generalization of the concept of Lie algebra.

The development of the theory of Lie algebras has led to numerous generalizations of the concept of Lie algebra. One such generalization arose in connection with the extension of the concept of Lie algebra to algebraic structures with $n$-ary multiplication laws. This generalization is called $n$-Lie algebra and it was proposed and developed by Filippov \cite{Filippov(1985)}. Independently of him, Nambu proposed a generalization of Hamiltonian mechanics based on the notion of an $n$-ary Poisson bracket \cite{Nambu(1973)}. It was later shown that the $n$-ary Poisson bracket in Nambu's generalization of Hamiltonian mechanics satisfies an $n$-Lie algebra identity (now called the Filippov-Jacobi or Fundamental Identity) and thus an $n$-ary Poisson bracket induces an $n$-Lie algebra structure on a vector space of smooth functions. It should be noted that the concept of $n$-Lie algebra turned out to be fruitful, and in the early 2000s this structure was used in the theory of M2-branes \cite{Bagger-Lambert(2007),Bagger-Lambert(2008)}. It is interesting to note that the quark model served as the motivation for Nambu to construct a generalization of Hamiltonian mechanics.

In this paper we propose a new approach to extending the concept of Lie algebra to algebraic structures with ternary multiplication laws. Our approach differs from the Filippov-Nambu approach, and to explain the difference, we briefly recall the main properties of $n$-Lie algebras. First, the $n$-ary Lie bracket of an $n$-Lie algebra is completely skew-symmetric, and second, the Filippov-Jacobi identity is an extension of the Leibniz rule to a double $n$-ary Lie bracket. The main examples of $n$-ary Lie brackets in $n$-Lie algebras are constructed using determinants. By this we mean that the theory of $n$-Lie algebras lacks an important construction that makes it possible to construct a Lie algebra using a {\em commutator}. It is well known that if $\cal A$ is an associative algebra over a field of real or complex numbers with multiplication $(u,v)\in {\cal A}\to u\cdot v\in{\cal A}$ then one can construct a Lie algebra by equipping $\cal A$ with Lie bracket defined by means of the commutator $[u,v]=u\cdot v-v\cdot u$. The commutator satisfies the Jacobi identity
\begin{equation}
\big[[u,v],w\big]+\big[[v,w],u\big]+\big[[w,u],v\big]=0,
\label{Jacobi identity}
\end{equation}
and $\cal A$ becomes the Lie algebra. This construction is very important, since it opens up the possibility of constructing a wide and important class of matrix Lie algebras.

Let us consider in more detail the above construction of a Lie algebra by means of the commutator. First of all, we are interested in the question of why the commutator satisfies the Jacobi identity.
It is easy to verify that if we expand all the double commutators on the left-hand side of the Jacobi identity (\ref{Jacobi identity}), then each permutation of three elements $u,v,w$ will appear in the resulting expression twice, once in the form of product $(u\cdot v)\cdot w$, the second time in the form of product
$u\cdot (v\cdot w)$, and these two products will have opposite signs. Thus, due to associativity $(u\cdot v)\cdot w=u\cdot (v\cdot w)$, we get zero. Our goal in this paper is to extend this construction to ternary algebras, that is, to construct a {\em ternary commutator} and find an {\em identity} for this ternary commutator (analogous to the Jacobi identity), which will be based on ternary associativity.

Let $\cal T$ be a vector space over the field of complex numbers $\mathbb C$ endowed with a ternary multiplication
$$
(a,b,c)\in {\cal T}\times{\cal T}\times{\cal T}\mapsto a\cdot b\cdot c\in {\cal T}.
$$
By ternary multiplication we mean a mapping that assigns to each triple of elements of complex vector space $\cal T$ the uniquely defined element of the same space $\cal T$, and this mapping is trilinear. In the case of ternary multiplication there are {\em two kinds of associativity}. A ternary multiplication is said to be associative of the first kind if it satisfies
\begin{equation}
(a\cdot b\cdot c)\cdot d\cdot f=a\cdot (b\cdot c\cdot d)\cdot f=a\cdot b\cdot (c\cdot d\cdot f),
\label{associativity of I kind}
\end{equation}
and associative of the second kind if
\begin{equation}
(a\cdot b\cdot c)\cdot d\cdot f=a\cdot (d\cdot c\cdot b)\cdot f=a\cdot b\cdot (c\cdot d\cdot f).
\label{associativity of II kind}
\end{equation}
Note that in the case of associativity of the first kind, shifting the round brackets from left to right does not change the order of the factors in a product, and in the case of associativity of the second kind, shifting the round brackets from left to right swaps elements $b$ and $d$. If we do not assume a vector space structure on $\cal T$, then a ternary multiplication on $\cal T$, satisfying associativity of the second kind (sometimes called generalized associativity), determines the structure of a {\em semi-heap} \cite{Wagner(1953)}. In what follows, a vector space $\cal T$ equipped with a ternary associative multiplication of the first or second kind will be referred to as a {\em ternary algebra}.

In the case of algebra with binary law of multiplication there are {\em two ways} of placing brackets in a product $u\cdot v\cdot w$ (which show in what order the elements $u,v,w$ are multiplied), that is, $(u\cdot v)\cdot w$ and $u\cdot (v\cdot w)$. If multiplication is associative then these two products are equal. Therefore, multiplying one of these two products by -1 and adding them together, we, due to associativity, get zero. This is the basis of the Jacobi identity.

Now in the case of ternary multiplication, we have {\em three ways} of placing round brackets in a product of five elements, that is,
\begin{equation}
(a\cdot b\cdot c)\cdot d\cdot f,\;\;\;\;a\cdot (b\cdot c\cdot d)\cdot f,\;\;\;\;
       a\cdot b\cdot (c\cdot d\cdot f),
\label{three products I}
\end{equation}
in the case of the associativity of the first kind and
\begin{equation}
(a\cdot b\cdot c)\cdot d\cdot f,\;\;\;\;a\cdot (d\cdot c\cdot b)\cdot f,\;\;\;\;
       a\cdot b\cdot (c\cdot d\cdot f),
\label{three products II}
\end{equation}
in the case of the associativity of the second kind. The main idea of the present paper is that we can find an analogue of the Jacobi identity for ternary algebras if we follow the above scheme, that is, we multiply each of the three products in (\ref{three products I}) or in (\ref{three products II}) by some number, add three resulting products and, by virtue of ternary associativity of the first or second kind, obtain zero. We think that using -1 in the ternary case looks unnatural, as it makes the whole construction asymmetrical and awkward. But the {\em cube roots of unity} fit perfectly into this scheme. Indeed there are {\em three} cube roots of unity $1,\omega,\overline\omega$, where $\omega$ is a primitive cube root of unity and $\overline\omega$ is its complex conjugate. Now if we multiply each product in (\ref{three products I}) or in (\ref{three products II}) by a cube root of unity (each by its own, different from the others) and add the resulting products, then by virtue of ternary associativity and the property of cube roots of unity $1+\omega+\overline\omega=0$ we get zero. For example
\begin{equation}
1\;(a\cdot b\cdot c)\cdot d\cdot f+\omega\;a\cdot (b\cdot c\cdot d)\cdot f+\overline\omega\;
       a\cdot b\cdot (c\cdot d\cdot f)=0,
\label{sum of products}
\end{equation}
and this equality can be considered as a ternary analog of the binary one $(u\cdot v)\cdot w-u\cdot (v\cdot w)=0.$

The above reasoning leads to the conclusion that in the case of ternary algebra $\cal T$ an analog of the binary commutator can be constructed by means of the cube roots of unity. We propose a ternary commutator which is a linear combination of all six permutations of its arguments and the coefficients of this linear combination are the cube roots of unity. Thus we endow a ternary algebra $\cal T$ with the ternary commutator defined by the following formula
\begin{eqnarray}
[a,b,c]=a\cdot b\cdot c + \omega\;b\cdot c\cdot a + \overline\omega\;c\cdot a\cdot b
        +c\cdot b\cdot a + \overline\omega\;b\cdot a\cdot c + \omega\; a\cdot c\cdot b.
        \label{introduction ternary commutator}
\end{eqnarray}
The structure of this ternary commutator is in line with the ideas, methods and structures developed in the papers \cite{Abramov-Kerner-LeRoy}, \cite{Abramov-Kerner-Liivapuu-Shitov}, \cite{Abramov-Liivapuu}, \cite{Kerner(2019)}.
The ternary commutator (\ref{introduction ternary commutator}) proposed in the present paper differs in its properties from a 3-Lie bracket of 3-Lie algebra. Indeed, the ternary commutator (\ref{introduction ternary commutator}) is not skew-symmetric and, therefore, the presence of two equal arguments does not make it identically zero. However, in the case where all three arguments are equal, it is identically zero. Here we see an analogy with the ternary generalization of the Pauli exclusion principle proposed by Kerner \cite{Kerner(2017)}. According to this principle, a wave function of a quantum system of quarks does not vanish in the case of two quarks with identical quantum characteristics, but it vanishes identically when the system contains three such quarks.

We find an identity for the ternary commutator (\ref{introduction ternary commutator}). It is natural to assume that, just as in the case of the Jacobi identity, which is based on the subgroup of cyclic permutations of three elements ${\mathbb Z}_3\subset S_3$, an identity for the ternary commutator (\ref{introduction ternary commutator}) should be based on a subgroup of the symmetric group $S_5$ and it really is. The identity we found is based on the general affine group $GA(1,5)\subset S_5$. Thus, the left side of identity has 20 terms, but in some cases of ternary multiplication with the commutativity with respect to first two arguments it decreases to 10 terms. The identity has the form
\begin{equation}
\circlearrowleft \Big(\big[[a,b,c],d,f\big]+\big[[a,d,b],f,c\big]
+\big[[a,f,d],c,b\big]+\big[[a,c,f],b,d\big]\Big) =0,\label{introduction_identity}
\end{equation}
where the symbol $\circlearrowleft$ stands for cyclic permutations of five elements and $a,b,c,d,f\in\cal T$.

Motivated by this result result we propose a notion of {\em ternary Lie algebra at cubic root of unity} or, more briefly, {\em ternary $\omega$-Lie algebra}, where $\omega$ is a primitive cube root of unity. We give definition of ternary $\omega$-Lie algebra and study its structure constants. The structure constants of a ternary $\omega$-Lie algebra is a (1,3)-tensor and we derive the system of equations for this tensor from the identity (\ref{introduction_identity}). We show that in the case of ternary $\omega$-Lie algebra with three generators the structure constants of this algebra are related to irreducible tensor representation of weight two of the rotation group.  Then we give several examples of ternary $\omega$-Lie algebras constructed using ternary algebras of rectangular and three-dimensional matrices. Thus, the concept of ternary $\omega$-Lie algebra proposed in this paper allows us to extend the theory of Lie algebras from square matrices to rectangular and three-dimensional matrices.

\section{Ternary Commutator and Its Symmetries}
In this section we explain  why we call the expression on the right-hand side of (\ref{introduction ternary commutator}) a ternary commutator. In addition, we describe the symmetries of the ternary commutator, define its conjugate ternary commutator, and derive a formula for the ternary commutator using sixth order roots of unity.

The concept of a commutator is closely related to a concept of commutativity. In the case of binary multiplication $u\cdot v$, two elements $u,v$ are called commuting if the equality $u\cdot v=v\cdot u$ holds. Now we introduce a commutator as an expression that vanishes on commuting elements, that is, $[u,v]=u\cdot v-v\cdot u$. In the case of $n$-ary multiplication, where $n>2$, commutativity can be defined in different ways, depending on how we interpret commutativity in the binary case.
For our purposes it is convenient to interpret the binary commutativity as follows: we split a product into two parts, and rearranging these parts does not change the value of the product. In this form, commutativity can be extended to $n$-ary multiplication laws. Assume we have an $n$-ary product $a_1\cdot a_2\cdot\ldots \cdot a_n$ and we split it into two parts
\begin{equation}
\underbrace{a_1\cdot a_2\cdot\ldots\cdot a_i}\cdot \underbrace{a_{i+1}\cdot a_{i+2}\cdot\ldots\cdot a_n},
\label{partition}
\end{equation}
where $i=1,2,\ldots,n-1$. We call $n$-ary multiplication {\em commutative} if for any $n$ elements $a_1,a_2,\ldots,a_n$ and for any partition of their product  into two parts (\ref{partition}) rearrangement of these two parts does not change the value of the product, i.e.
$$
a_1\cdot a_2\cdot\ldots\cdot a_i\cdot a_{i+1}\cdot a_{i+2}\cdot\ldots\cdot a_n=
   a_{i+1}\cdot a_{i+2}\cdot\ldots\cdot a_n\cdot a_1\cdot a_2\cdot\ldots\cdot a_i.
$$
Applying this notion of commutativity to the special case of a ternary multiplication, i.e. assuming that a ternary multiplication is commutative in the sense just defined we will have the relations
\begin{eqnarray}
&& \underbrace{a}\cdot\underbrace{b\cdot c}=\underbrace{b\cdot c}\cdot \underbrace{a},\;\;
   \underbrace{a\cdot b}\cdot\underbrace{c}=\underbrace{c}\cdot \underbrace{a\cdot b},\label{commutativity 1}\\
&& \underbrace{c}\cdot\underbrace{b\cdot a}=\underbrace{b\cdot a}\cdot \underbrace{c},\;\;
   \underbrace{c\cdot b}\cdot\underbrace{a}=\underbrace{a}\cdot \underbrace{c\cdot b}.\label{commutativity 2}
\end{eqnarray}
Thus, these relations show that the notion of commutativity in the special case of a ternary multiplication is equivalent to the fact that any cyclic permutation of arguments in a triple product does not change the value of this product. We have two sets of equal products (cyclic permutations of three elements), and each set contains {\em three products}. Now, in order to construct a ternary commutator, we should combine all six products into a linear combination in such a way that, first, we should use the cube roots of unity, and, second, when the commutativity conditions (\ref{commutativity 1}), (\ref{commutativity 2}) are satisfied, this combination would
vanish. Hence we define the ternary commutator as follows
\begin{equation}
a\cdot b\cdot c + \omega\;b\cdot c\cdot a + \overline\omega\;c\cdot a\cdot b
        +c\cdot b\cdot a + \overline\omega\;b\cdot a\cdot c + \omega\; a\cdot c\cdot b,
\end{equation}
where $a,b,c$ are elements of a ternary algebra and $\omega$ is a primitive 3rd order root of unity. In what follows we will call this expression {\em ternary commutator} and denote it using square brackets, that is,
\begin{equation}
[a,b,c]=a\cdot b\cdot c + \omega\;b\cdot c\cdot a + \overline\omega\;c\cdot a\cdot b
        +c\cdot b\cdot a + \overline\omega\;b\cdot a\cdot c + \omega\; a\cdot c\cdot b.
\label{ternary commutator}
\end{equation}
The ternary commutator can be constructed by analogy with the binary one using geometric reasoning. In the case of binary multiplication $u\cdot v$, it is natural to place the factors $u,v$ at the ends of the segment and form their product by taking as the first factor the element that stands at the left end of the segment. The permutation of the factors in a product $u\cdot v$ corresponds to the rotation of the segment around its center by an angle of $\pi$. Thus, the binary commutator can be interpreted in such a way that we take the product $u\cdot v$ determined by the initial position of the segment and add to it the product determined by the segment rotated by the angle of $\pi$, multiplied by the coefficient $e^{i\pi}$, that is, $u\cdot v+e^{i\pi}\;v\cdot u=u\cdot v-v\cdot u=[u,v]$.

In the case of ternary multiplication, we should use a regular triangle to graphically represent a triple product $a\cdot b\cdot c$. We will arrange the three factors $a,b,c$ of this product at the vertices of a triangle, placing the first factor $a$ at the lower right vertex and going around the triangle clockwise. Then rotating the triangle around its center by an angle of $2\pi/3$ counterclockwise will give us the first cyclic permutation $b\cdot c\cdot a$, and rotating it by an angle of $4\pi/3$ will give us the second $c\cdot a\cdot b$. So they must enter into the expression for the ternary commutator with the factors $e^{2\pi i/3}=\omega,e^{4\pi i/3}=\overline\omega$. The second part of the expression for the ternary commutator is obtained by mirroring the described construction.
\vskip.3cm
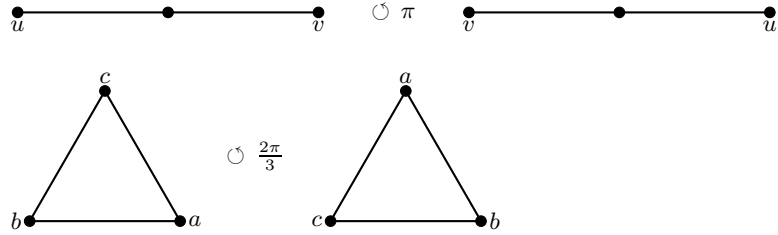
\begin{tikzpicture}
    \draw[thick] (-4,0) -- (0,0);

    \filldraw[black] (-4,0) circle (2pt); 
    \filldraw[black] (0,0) circle (2pt);  
    \filldraw[black] (-2,0) circle (2pt); 

    \node[below] at (-4,0) {$u$}; 
    \node[below] at (0,0) {$v$};  

    \draw[thick] (2,0) -- (6,0);

    \filldraw[black] (2,0) circle (2pt); 
    \filldraw[black] (6,0) circle (2pt);  
    \filldraw[black] (4,0) circle (2pt);  

    \node[below] at (2,0) {$v$}; 
    \node[below] at (6,0) {$u$}; 

    \node at (1,0) {\(\circlearrowleft\,\pi\)};
\end{tikzpicture}
\vskip.3cm

\begin{tikzpicture}

\draw[thick] (0,0) -- (2,0) -- (1,1.732) -- cycle;
\filldraw[black] (0,0) circle (2pt) node[left] {$b$};
\filldraw[black] (2,0) circle (2pt) node[right] {$a$};
\filldraw[black] (1,1.732) circle (2pt) node[above] {$c$};

\draw[thick] (4,0) -- (6,0) -- (5,1.732) -- cycle;
\filldraw[black] (4,0) circle (2pt) node[left] {$c$};
\filldraw[black] (6,0) circle (2pt) node[right] {$b$};
\filldraw[black] (5,1.732) circle (2pt) node[above] {$a$};

\node at (3,0.866) { $\quad\circlearrowleft \, \frac{2\pi}{3}\quad$};
\end{tikzpicture}

\noindent
It is useful to define the conjugate ternary commutator as follows
\begin{eqnarray}
[a,b,c]^\ast=a\cdot b\cdot c + \overline\omega\;b\cdot c\cdot a + \omega\;c\cdot a\cdot b
        +c\cdot b\cdot a + \omega\;b\cdot a\cdot c + \overline\omega\; a\cdot c\cdot b.
        \label{ternary conjugate commutator}
\end{eqnarray}
Then we have
$$
[a,b,c]^\ast = [c,b,a],\;\;\;([a,b,c]^\ast)^\ast=[a,b,c].
$$
It is easy to verify by direct calculation that the ternary commutator (\ref{ternary commutator}) and its conjugate (\ref{ternary conjugate commutator}) transform under cyclic permutations of their arguments as follows
\begin{eqnarray}
&& [a,b,c]=\omega\;[b,c,a],\;\;[a,b,c]=\overline\omega\;[c,a,b],\label{cyclic}\\
&& [a,b,c]^\ast=\overline\omega\,{[b,c,a]^\ast},\; [a,b,c]^\ast=\omega\,{[c,a,b]^\ast}.\label{noncyclic}
\end{eqnarray}
From (\ref{cyclic}) it follows that the sum of three ternary commutators obtained by cyclic permutations of arguments is equal to zero, that is,
\begin{equation}
[a,b,c]+[b,c,a]+[c,a,b]=0.
\label{sum of cyclic permutations is zero}
\end{equation}
Concerning this important property of the ternary commutator (\ref{ternary commutator}) we have to make three remarks. The first remark concerns Lie triple systems that arose in differential geometry in connection with the study of totally geodesic submanifolds. Although the property (\ref{sum of cyclic permutations is zero}) of the ternary commutator (\ref{ternary commutator}) has the same form as one of the requirements in the definition of a Lie triple system, the ternary commutator (\ref{ternary commutator}) is not a Lie triple system because it is not skew-symmetric in the first two arguments and does not satisfy the Filippov-Jacobi identity.

The second remark concerns 3-Lie algebras. Ternary Lie bracket in 3-Lie algebra is totally skew-symmetric in its arguments and thus in general it does not satisfy the equation (\ref{sum of cyclic permutations is zero}). In the next section we will find an identity for the ternary commutator (\ref{ternary commutator}) and see that this identity is different from the Filippov-Jacobi identity. Thus, the ternary commutator (\ref{ternary commutator}) induces on a ternary algebra $\cal T$ a structure different from a 3-Lie algebra.

The third remark concerns a relation with theoretical physics. It is easy to see that two equal arguments in our ternary commutator (\ref{ternary commutator}) does not make it vanish identically. But in the case where all three arguments $a,b,c$ are equal, our ternary commutator vanishes identically, i.e. $[a,a,a]=0$. Here we see a possible connection with the ternary generalization of the Pauli exclusion principle proposed by Kerner \cite{Kerner(2017)}. We will discuss this connection in the Discussion section.

Recall that Nambu in the paper \cite{Nambu(1973)} devoted to the generalization of Hamiltonian mechanics considered the skew-symmetric ternary commutator
$$
[A,B,C]=A\,B\,C+B\,C\,A+C\,A\,B-C\,B\,A-B\,A\,C-A\,C\,B,
$$
where $A,B,C$ are elements of some associative (binary) algebra. This version of a ternary commutator can be considered as a direct extension of the skew-symmetry of the binary commutator to the case of ternary multiplication. However, to our knowledge, no analogue of the Jacobi identity based on ternary associativity has been found for such a ternary commutator. It is interesting that our ternary commutator (\ref{ternary commutator}) can also be written in a form where the three even permutations have a plus sign, and the three odd permutations have a minus sign. For this purpose we will need a primitive sixth root of unity which will be denoted by $\varepsilon$. We take $\omega=\varepsilon^2, \overline\omega=\varepsilon^4$. Among other relations we will have
\begin{eqnarray}
    \varepsilon+\overline\varepsilon=1,\;\;\omega=-\overline\varepsilon,\;\;
             \overline\omega=-\varepsilon.
\end{eqnarray}
Now we can write the ternary commutator (\ref{ternary commutator}) in the following form
\begin{eqnarray}
[a,b,c] =a\cdot b\cdot c-\varepsilon\;b\cdot a\cdot c+\varepsilon^2\;b\cdot c\cdot a
                   -\varepsilon^3\; c\cdot b\cdot a + \varepsilon^4\;c\cdot a\cdot b
                   -\varepsilon^5\;a\cdot c\cdot b.\label{ternary commutator with epsilon}
\end{eqnarray}
In this formula, even permutations have a plus sign and are multiplied by even powers of the sixth root of unity $\varepsilon$, and odd permutations have a minus sign and are multiplied by odd powers of the sixth root of unity. Now the symmetries of the ternary commutator can be written in the form
\begin{eqnarray}
&& [a,b,c] = \varepsilon^2\,[b,c,a]=\varepsilon^4\,[c,a,b],\;\;\;\;
                [a,b,c]^\ast = \varepsilon^4\,[b,c,a]^\ast=\varepsilon^2\,[c,a,b]^\ast,\nonumber\\
&& [a,b,c]=-\varepsilon\;[b,a,c]^\ast,\;\;[a,b,c] =-\varepsilon^3\; [c,b,a]^\ast,\;\;
      [a,b,c]=-\varepsilon^5\;[a,c,b]^\ast.\nonumber
\end{eqnarray}
The formula (\ref{ternary commutator with epsilon}) can be written in the form
\begin{equation}
[a,b,c] =(a\cdot b\cdot c-\varepsilon\;b\cdot a\cdot c)+\omega\;(b\cdot c\cdot a
                   -\varepsilon\; c\cdot b\cdot a) + \overline\omega\;(c\cdot a\cdot b
                   -\varepsilon\;a\cdot c\cdot b).
                   \label{ternary commutator form 3}
\end{equation}
We can use this formula to justify the term "ternary commutator" that we use in relation to expression on the right-hand side of (\ref{ternary commutator}). In the above formula, each of the three terms enclosed in round brackets can be interpreted as measuring the non-commutativity of the ternary multiplication with respect to the first two arguments in relation to the last which does not change its position. Geometrically, it would be convenient to depict the three elements $a,b,c$ of the ternary algebra $\cal T$ as the vertices of a regular triangle. Then the above formula "measures" the non-commutativity of a ternary multiplication on each side of the triangle {\em with respect to the opposite vertex}. Thus, geometrically, the transition from binary multiplication, where two factors can be represented as points on a line, to ternary multiplication can be described as we leave a line and go to a plane, figuratively speaking. This explains why the above formula contains sixth roots of unity and conjugation. To measure ternary non-commutativity correctly, we need to use plane rotations and reflections. It should be noted here that in \cite{Zapata-Arsiwalla-Beynon(2024)} the authors develop an interesting graphical and diagrammatic approach for representing ternary associative multiplication using triangles in the plane.

In particular, if a ternary multiplication is commutative with respect to some pair of its arguments, for example the first pair, that is, $a\cdot b\cdot c=b\cdot a\cdot c$, then formula (\ref{ternary commutator form 3}) reduces to a shorter form containing only three terms. Indeed, we have
\begin{eqnarray}
&& [a,b,c]=a\cdot b\cdot c+\omega\,b\cdot c\cdot a+\overline\omega\,c\cdot a\cdot b+
            c\cdot b\cdot a+\overline\omega\,b\cdot a\cdot c+\omega\,a\cdot c\cdot b\nonumber\\
&&\;\;\;\;\qquad = (1+\overline\omega)\,a\cdot b\cdot c+(1+\omega)\,b\cdot c\cdot a+(\omega+\overline\omega)\,c\cdot a\cdot b\nonumber\\
&&\;\quad\qquad = -\omega\,(a\cdot b\cdot c+\omega\,b\cdot c\cdot a+\overline\omega\,c\cdot a\cdot b).
\end{eqnarray}

\section{General Affine Group, Basic Identity and\\ Ternary Lie Algebra at Cube Root of Unity}
The concept of a Lie algebra consists of two important components, where the first is a Lie bracket (or, in particular, the binary commutator) with its properties with respect to permutations of arguments and the Jacobi identity. Since we have the ternary commutator defined and considered in the previous section, our goal now is to find an identity for the ternary commutator (\ref{ternary commutator}), based on ternary associativity. Following the analogy with the binary commutator, we could estimate how many terms a possible identity could contain. If we consider the binary case then each double commutator, when expanded, yields four products. But if we expand all the double commutators at the left-hand side of identity, then in the resulting expression each product of three elements (totally we have six permutations) will occur twice (the brackets are either on the left or on the right). Thus, we will have twelve products on the left-hand side of identity. Dividing twelve by four we conclude that an identity consists of three double commutators and this is so in the case of the Jacobi identity.

A similar calculation can be made in the case of the ternary commutator (\ref{ternary commutator}). If we expand the double ternary commutator $[[a,b,c],d,f]$ we get thirty-six terms. On the other hand, we have one hundred and twenty permutations of five elements. Due to ternary associativity, each permutation must occur at least three times (brackets on the left, in the center and on the right) with coefficients $1,\omega,\overline\omega$. Thus, dividing three hundred and sixty by thirty-six gives ten. Note that this is the minimum number of terms in a possible identity. Also note that in this calculation we have not taken into account such an important structure of the ternary commutator as conjugation. Obviously, if we take this structure into account, we will have to double the number of terms in the identity, i.e. we can expect that a possible identity will contain twenty terms.

Since an identity we are looking for is a sum of double ternary commutators of the form $[[a,b,c],d,f]$, the second assumption, which seems very natural, is that an identity must be based on a subgroup of symmetric group $S_5$. Taking into account the above, we come to the conclusion that there are two potential candidates for the subgroups of the symmetric group $S_5$, these are the dihedral group $D_{10}$ (10 elements) or the general affine group $GA(1,5)$ (20 elements). Moreover, the dihedral group is a subgroup of the general affine group, that is, $D_{10}\subset GA(1,5)$.

The general affine group $GA(1,5)$ has several different representations. In this article we will use the representation of this group by permutations of five elements. The minimal set of permutations that generates the entire group consists of two cycles, which we denote as follows
$$
\sigma=(1\; 2\; 3\; 4\; 5),\;\;\tau=(2\;4\;5\;3).
$$
Hence
\begin{equation}
\sigma(1)=2,\;\sigma(2)=3,\;\sigma(3)=4,\;\sigma(5)=1,
\end{equation}
and
\begin{equation}
\tau(1)=1,\;\tau(2)=4,\;\tau(3)=2,\;\tau(4)=5,\;\tau(5)=3.
\end{equation}
Then
$$
GA(1,5)=<\sigma,\tau\;|\; \sigma^5=e, \tau^4=e, \tau\,\sigma\,\tau^{-1}=\sigma^2>,
$$
where $e$ is the identity element of the group $GA(1,5).$ All elements of the group can be written in the following form
\begin{eqnarray}
&& e,\; \sigma,\; \sigma^2,\;\sigma^3,\;\sigma^4,\label{first 5}\\
&& \tau,\;\tau\sigma^3,\;\tau\sigma,\;\tau\sigma^4,\;\tau\sigma^2,\label{second 5}\\
&& \tau^2,\;\tau^2\sigma^4,\;\tau^2\sigma^3,\;\tau^2\sigma^2,\;\tau^2\sigma,\label{third 5}\\
&& \tau^3,\;\tau^3\sigma^2,\;\tau^3\sigma^4,\;\tau^3\sigma,\;\tau^3\sigma^3.\label{fourth 5}
\end{eqnarray}
We will use this representation to write the identity. In this representation, all elements of the general affine group are divided into four sets (\ref{first 5}), (\ref{second 5}), (\ref{third 5}), (\ref{fourth 5}) and in each of these sets the second element is obtained by a cyclic permutation of five elements in the first, the third by a cyclic permutation in the second, and so on. For a more compact representation of an identity, we will use the symbol $\circlearrowleft$. This symbol means that an expression that follows contains five elements and must be subjected to the following procedure. One should form the sum of five expressions, starting with the initial one and where each subsequent one is a cyclic permutation of five elements of the previous one. Thus
\begin{eqnarray}
\circlearrowleft \big[[a,b,c],d,f\big]\!\!\! &=&\!\!\!
     \big[[a,b,c],d,f\big]+\big[[b,c,d],f,a\big]+\big[[c,d,f],a,b\big]+\big[[d,f,a],b,c\big]\nonumber\\
     &&
     \qquad\qquad\qquad\qquad\qquad\qquad\qquad\qquad\qquad\qquad\qquad+\big[[f,a,b],c,d\big], \nonumber
\end{eqnarray}
where $a,b,c,d,f$ are elements of a ternary algebra $\cal T$.
\begin{theorem}
Let $\cal T$ be a ternary algebra. Then for any $a,b,c,d,f\in\cal T$ the ternary commutator (\ref{ternary commutator}) and its conjugate (\ref{ternary conjugate commutator}) have the property
$$
[a,b,c]=\omega\,[b,c,a]=\overline\omega\,[c,a,b],\;\;\;
  [a,b,c]^\ast=\overline\omega\,{[b,c,a]^\ast}=\omega\,{[c,a,b]^\ast}
$$
and the ternary commutator satisfies the identity
\begin{eqnarray}
\circlearrowleft \Big(\big[[a,b,c],d,f\big]+\big[[a,d,b],f,c\big]
+\big[[a,f,d],c,b\big]+\big[[a,c,f],b,d\big]\Big) =0.\nonumber
\end{eqnarray}
\label{theorem identity}
\end{theorem}
In what follows we will call the above identity that is the statement of Theorem \ref{theorem identity} the {\em basic identity}. We can prove Theorem \ref{theorem identity} by direct computation, that is, by applying formula (\ref{ternary commutator}) twice to each term of the basic identity and using a ternary associativity of multiplication. We carried out this computation using a computer program containing a non-commutative symbolic calculus. The computer program we use makes it possible to study the structure of the basic identity. A study of the structure of the basic identity shows that it holds due to the reasoning based on the formulas (\ref{three products I}), (\ref{three products II}), (\ref{sum of products}). Let us denote $a=a_1,b=a_2,c=a_3,d=a_4,f=a_5$. The computer program allows us to find in which terms of the basic identity a particular product of elements $a_1,a_2,a_3,a_4,a_5$ appears, with what coefficient and how the round brackets are placed in it. For example the product $a_1\cdot a_2\cdot a_3\cdot a_4\cdot a_5$ appears six times as follows
\begin{eqnarray}
&& [[a_1,a_2,a_3],a_4,a_5],\;\;\;\;\;\;\;\;\;[[a_2,a_3,a_4],a_5,a_1],\;\;\;\;\;\;\;\;\;\;\;\;[[a_3,a_4,a_5],a_1,a_2],\nonumber\\
&&\!\!\! (a_1\cdot a_2\cdot a_3)\cdot a_4\cdot a_5,\;\;\;\;\;\;
     \overline\omega\; a_1\cdot (a_2\cdot a_3\cdot a_4)\cdot a_5, \;\;\;\;\;\;
        \omega\;a_1\cdot a_2\cdot (a_3\cdot a_4\cdot a_5),\nonumber
\end{eqnarray}
and
\begin{eqnarray}
&& [[a_5,a_4,a_3],a_2,a_1],\;\;\;\;\;\;\;\;\;[[a_4,a_3,a_2],a_1,a_5],\;\;\;\;\;\;\;\;\;\;\;\;[[a_3,a_2,a_1],a_5,a_4],\nonumber\\
&&\!\!\! a_1\cdot a_2\cdot (a_3\cdot a_4\cdot a_5),\;\;\;\;\;\;
     \overline\omega\; a_1\cdot (a_2\cdot a_3\cdot a_4)\cdot a_5, \;\;\;\;\;\;
        \omega\;(a_1\cdot a_2\cdot a_3)\cdot a_4\cdot a_5.\nonumber
\end{eqnarray}
Here in the first line we show the double ternary commutators of the basic identity and below them we show in which form, that is, the coefficient and position of round brackets, the product $a_1\cdot a_2\cdot a_3\cdot a_4\cdot a_5$ appears in the corresponding double commutator. Adding up the six terms obtained in this case and assuming associativity of the first kind, we get zero in total.

In the case of ternary associativity of the second kind, in addition to the above table we should consider the set of those double ternary commutators on the left-hand side of the basic identity that contain the product $a_1\cdot a_4\cdot a_3\cdot a_2\cdot a_5$. They are summarized in the following table
\begin{eqnarray}
&& [[a_3,a_1,a_4],a_2,a_5],\;\;\;\;\;\;\;\;\;[[a_2,a_3,a_4],a_5,a_1],\;\;\;\;\;\;\;\;\;\;\;\;[[a_2,a_5,a_3],a_1,a_4],\nonumber\\
&&\!\!\! \omega\;(a_1\cdot a_4\cdot a_3)\cdot a_2\cdot a_5,\;\;\;\;
     \overline\omega\;a_1\cdot (a_4\cdot a_3\cdot a_2)\cdot a_5, \;\;\;\;\;\;
        a_1\cdot a_4\cdot (a_3\cdot a_2\cdot a_5),\nonumber
\end{eqnarray}
and
\begin{eqnarray}
&& [[a_4,a_1,a_3],a_5,a_2],\;\;\;\;\;\;\;\;\;[[a_4,a_3,a_2],a_1,a_5],\;\;\;\;\;\;\;\;\;\;\;\;[[a_3,a_5,a_2],a_4,a_1],\nonumber\\
&&\!\!\! (a_1\cdot a_4\cdot a_3)\cdot a_2\cdot a_5,\;\;\;\;\;\;
     \overline\omega\; a_1\cdot (a_4\cdot a_3\cdot a_2)\cdot a_5, \;\;\;\;\;\;
        \omega\;a_1\cdot a_4\cdot (a_3\cdot a_2\cdot a_5).\nonumber
\end{eqnarray}
A comparison of the columns in the center of these tables immediately shows that in the case of ternary associativity of the second kind we obtain the same type of sum, which is equal to zero.

Thus the basic identity consists of 20 double ternary commutators. The general affine group, considered as a subgroup of the permutations of the symmetric group 1, is generated by two cycles. The general affine group, considered as a subgroup of the permutations of the symmetric group $S_5$, is generated by two cycles $\sigma,\tau$. The double ternary commutators
\begin{equation}
\big[[a_1,a_2,a_3],a_4,a_5\big],\;\big[[a_1,a_4,a_2],a_5,a_3\big],\;\big[[a_1,a_5,a_4],a_3,a_2\big],\;\big[[a_1,a_3,a_5],a_2,a_4\big]
\label{Z_4}
\end{equation}
are determined by the permutations $e,\tau,\tau^2,\tau^3$, that is, by the first elements of the general affine group $GA(1,5)$ in the formulas (\ref{first 5}),(\ref{second 5}),(\ref{third 5}),(\ref{fourth 5}). The cyclic permutations of the double ternary commutators (\ref{Z_4}) are determined by the elements in (\ref{first 5}),(\ref{second 5}),(\ref{third 5}),(\ref{fourth 5}) respectively, starting from the second element. Note that the elements in (\ref{first 5}), (\ref{third 5}) form the dihedral subgroup $D_{10}$ of the symmetric group $S_5$.

Theorem \ref{theorem identity} provides a motivation for introduction of the following notion.
\begin{definition}
Let $\cal L$ be a vector space over the field of complex numbers. Then $\cal L$ is said to be a {ternary Lie algebra at cube root of unity}, where $\omega$ is a primitive cube root of unity, if $\cal L$ is endowed with a ternary bracket $(x,y,z)\in {\cal L}\times{\cal L}\times{\cal L}\mapsto [x,y,z]\in {\cal L}$ which transforms under the cyclic permutations of its arguments as follows
\begin{equation}
[x,y,z]=\omega\;[y,z,x]=\overline\omega\;[z,x,y],\;\;
     [x,y,z]^\ast=\overline\omega\;[y,z,x]^\ast=\omega\;[z,x,y]^\ast,
\label{omega-symmetry}
\end{equation}
where $[x,y,z]^\ast=[z,y,x]$, and satisfies the identity
\begin{eqnarray}
\circlearrowleft \Big(\big[[x,y,z],u,v\big]+\big[[x,u,y],v,z\big]
+\big[[x,v,u],z,y\big]+\big[[x,z,v],y,u\big]\Big) =0.\label{identity}
\end{eqnarray}
\end{definition}
In this paper, to simplify the terminology, a ternary Lie algebra at cube root of unity will be referred to as a ternary $\omega$-Lie algebra. The property (\ref{omega-symmetry}) will be referred to as $\omega$-symmetry of ternary bracket of a ternary $\omega$-Lie algebra. The identity (\ref{identity}) will be called as before the {\em basic identity}.

Let $\cal L$ be a ternary $\omega$-Lie algebra, where $\cal L$ is an $n$-dimensional vector space, and $e_1,e_2,\ldots,e_n$ be a basis for a vector space $\cal L$. In analogy with the binary case we introduce the structure constants of a ternary $\omega$-Lie algebra as follows
\begin{equation}
[e_i,e_k,e_l]=C_{ikl}^m\;e_m,\;\;\;[e_i,e_k,e_l]^\ast={\tilde C}_{ikl}^m\;e_m
\label{structure constants}
\end{equation}
where $C^m_{ikl},{\tilde C}^m_{ijk}$ will be referred to as structure constants of a ternary $\omega$-Lie algebra $\cal L$. It is easy to see that ${\tilde C}^m_{ijk}=C^m_{kji}$. In (\ref{structure constants}) we used the Einstein convention of summation over repeated indices. Obviously the structure constants of a ternary $\omega$-Lie algebra can be considered as a complex-valued tensor of type $(1,3)$. This tensor has the $\omega$-symmetry with respect to cyclic permutations of its three subscripts, that is,
\begin{equation}
C^m_{ikl}=\omega\;C^m_{kli}=\overline\omega\;C^m_{lik},\;\;\;
     {\tilde C}^m_{ikl}=\overline\omega\;{\tilde C}^m_{kli}=\omega\;{\tilde C}^m_{lik}.
\label{structure constants omega-symmetry}
\end{equation}
It follows that for every value of the superscript $m=1,2,\ldots,n$ the structure constants of a ternary $\omega$-Lie algebra $\cal L$, that is, both $C^m_{ijk}$ and $\tilde C^m_{ijk}$, satisfy the equation
\begin{equation}
T_{ijk}+T_{jki}+T_{kij}=0,
\label{tensors 1,3}
\end{equation}
where $T_{ijk}$ is a covariant tensor of order 3. It is evident that the third order covariant tensors defined on the vector space $\cal L$, which satisfy the equation (\ref{tensors 1,3}), form the subspace in the vector space of covariant tensors of order 3. This subspace will be denoted by $\tensors$.

The formulas (\ref{structure constants omega-symmetry}) clearly show that the for any superscript $m$ the structure constants $C^m_{ijk}$, $\tilde C^m_{ijk}$ are the eigenvectors of the linear operator in $\tensors$ induced by the cyclic permutation $(1\;2\; 3)$ with eigenvalues $\omega, \overline\omega$, respectively. Thus
\begin{equation}
\tensors=\1tensors\oplus\2tensors,
\label{decomposition}
\end{equation}
where
$$
\1tensors=\{T_{ijk}\in\tensors: T_{ijk}=\omega T_{jki}\},\;\;\;\;
    \2tensors=\{T_{ijk}\in\tensors: T_{ijk}=\overline\omega\; T_{jki}\}.
$$
Thus, for each value of the superscript $m$, the structure constants $C^m_{ijk}$ of a ternary $\omega$-Lie algebra $\cal L$ belong to subspace $\1tensors$, and the structure constants $\tilde C^m_{ijk}$ belong to subspace $\2tensors$.

Here we would like to note an important connection between the structure constants of a three-dimensional ternary Lie algebra and irreducible representations of the rotation group. Let $n=3$, i.e. we are considering a three-dimensional ternary Lie algebra. Let $A=(A^i_j)\in \mbox{SO}(3)$ be a real orthogonal matrix with determinant 1. Then the formula
\begin{equation}
T_{ijk}\;\to\; T^\prime_{prs}=A_p^i A_r^j A_s^k\; T_{ijk},
\label{representation}
\end{equation}
where $T_{ijk},T^\prime_{prs}\in\tensors$, defines a linear representation of the rotation group $ \mbox{SO}(3)$ in the space $\tensors$. If we add to equation (\ref{tensors 1,3}) the condition of tracelessness of a tensor $T_{ijk}$ for any pair of subscripts, formula (\ref{representation}) defines a twice repeated irreducible representation of the rotation group in the corresponding subspace of 3rd order covariant tensors \cite{Gelfand-Minlos-Shapiro}. Now the decomposition (\ref{decomposition}) splits this two-fold irreducible representation into two irreducible ones, respectively in subspaces $\1tensors$ and $\2tensors$ (with the additional condition that a tensor $T_{ijk}$ is traceless).
Note that the subspace of traceless tensors in $\1tensors$ is a five-dimensional Hermitian space and the explicit description of this space can be found in \cite{Abramov-Liivapuu}. In the next paper we plan to use this connection with irreducible representations of the rotation group to classify three-dimensional ternary $\omega$-Lie algebras.

It follows from the basic identity (\ref{identity}) that the structure constants of a ternary $\omega$-Lie algebra $\cal L$ satisfy the system of equations
\begin{equation}
\circlearrowleft (C^m_{\underline{i}\,\underline{k}\,\underline{l}}\,C^p_{m\,\underline{r}\,\underline{s}}+C^m_{\underline{i}\,\underline{r}\,\underline{k}}\,C^p_{m\,\underline{s}\,\underline{l}}+C^m_{\underline{i}\,\underline{s}\,\underline{r}}\,C^p_{m\,\underline{l}\,\underline{k}}+C^m_{\underline{i}\,\underline{l}\,\underline{s}}\,C^p_{m\,\underline{k}\,\underline{r}})=0.
\label{identity for structure constants}
\end{equation}
In this formula, the symbol $\circlearrowleft$ means that in an expression that follows it, we should perform the five cyclic permutations of the underlined subscripts and then take the sum of obtained expressions. For instance, if we apply $\circlearrowleft$ to the first term in (\ref{identity for structure constants}) then we get
$$
\circlearrowleft C^m_{\underline{i}\,\underline{k}\,\underline{l}}\,C^p_{m\,\underline{r}\,\underline{s}}=C^m_{ikl}\,C^p_{mrs}+C^m_{klr}\,C^p_{msi}+C^m_{lrs}\,C^p_{mik}+C^m_{rsi}\,C^p_{mkl}+C^m_{sik}\,C^p_{mlr}.
$$
We have the simplest case of a ternary $\omega$-Lie algebra in dimension 2.
It is easy to verify that if the vector space of a ternary $\omega$-Lie algebra has dimension 2, that is, the ternary $\omega$-Lie algebra has two generators $e_1,e_2$, the basic identity (\ref{identity}) does not impose any additional conditions, that is, it is satisfied due to the $\omega$-symmetries of ternary bracket. Using this, it is easy to show that in dimension 2, up to isomorphism, there is only one ternary $\omega$-Lie algebra, which is given by the following non-trivial commutation relations
\begin{equation}
[e_1,e_2,e_1]=e_2,\;\;\;\;[e_2,e_1,e_2]=e_1.
\end{equation}
We denote this 2-dimensional ternary $\omega$-Lie algebra by ${\cal L}_2$.
\section{Examples of Ternary Lie Algebra at Cube Root of Unity}
In this section we give some important examples of ternary associative algebras and consider ternary $\omega$-Lie algebras that are induced by the ternary commutator (\ref{ternary commutator}). A wide class of ternary associative algebras can be constructed using square matrices. Indeed, if $A,B,C$ are square matrices of order $n$, we can consider the ternary product $A\,B\,C$. This definition is correct, since matrix multiplication is associative. Obviously, in this case we obtain ternary multiplication with associativity of the first kind. However, from our point of view, this example is of little interest from the ternary point of view, because, firstly, ternary multiplication is constructed using binary, that is, binary is more fundamental than ternary, and, secondly, for square matrices there is a deeply developed theory of (binary) Lie algebras. Therefore, in this section we will consider examples of ternary algebras constructed using either rectangular (two-dimensional) matrices or cubic (three-dimensional) matrices. Thus, the notion of a ternary $\omega$-Lie algebra proposed in this paper can be considered as an extension of the concept of Lie algebra to rectangular and cubic matrices. Note that, firstly, the ternary multiplications considered in this section cannot be reduced to binary ones, and, secondly, they are associative of the second kind.

One of the simplest examples of ternary algebra with associativity of the second kind is an $n$-dimensional complex vector space $\mathbb C^n$ with a bilinear symmetric form $B$ defined on it. Then the ternary multiplication in $\mathbb C^n$ will be defined as follows
\begin{equation}
x\cdot y\cdot z=B(x,y)\,z.
\label{ternary product of vectors}
\end{equation}
It is easy to verify that this ternary product is associative of the second kind. Indeed
\begin{eqnarray}
&& (x\cdot y\cdot z)\cdot u\cdot v=(B(x,y)\,z)\cdot u\cdot v=
                 B(x,y)\,B(z,u)\,v,\nonumber\\
&& x\cdot (u\cdot z\cdot y)\cdot v=x\cdot (B(u,z)\,y)\cdot v=
                {B(u,z)}\,B(x,y)\,v,\nonumber\\
&& x\cdot y\cdot (z\cdot u\cdot v)=x\cdot y\cdot (B(z,u)\,v)=
                 B(z,u)\,B(x,y)\,v,\nonumber
\end{eqnarray}
and, due to the symmetry $B(z,u)={B(u,z)}$, we see that all three products are equal.
Hence if we endow a vector space $\mathbb C^n$ with the ternary commutator (\ref{ternary commutator}) then according to Theorem \ref{theorem identity} it becomes a ternary $\omega$-Lie algebra. The ternary commutator in this case can be written as follows
\begin{eqnarray}
[x,y,z] \!\!\!&=&\!\!\! x\cdot y\cdot z+\omega\,y\cdot z\cdot x+\overline\omega\,z\cdot x\cdot y+
           z\cdot y\cdot x+\overline\omega\,y\cdot x\cdot z+\omega\,x\cdot z\cdot y \nonumber\\
    \!\!\!&=&\!\!\! B(x,y)\, z+\omega\,B(y,z)\, x+\overline\omega\,B(z,x)\, y+
           B(z,y) x+\overline\omega\,B(y,x)\, z+\omega\,B(x,z)\, y\nonumber\\
    \!\!\!&=&\!\!\! (1+\overline\omega)\,B(x,y)\,z+(1+\omega)\,B(y,z)\,x+(\omega+\overline\omega)\,B(z,x)\,y\nonumber\\
    \!\!\!&=&\!\!\! -\big(B(z,x)\,y+\omega\,B(x,y)\,z+\overline\omega\,B(y,z)\,x\big).\nonumber
\end{eqnarray}
Omitting the irrelevant factor $-1$, we can consider the ternary commutator (\ref{ternary commutator}) and its conjugate in a reduced form
\begin{eqnarray}
&&[x,y,z]\;=B(z,x)\,y+\omega\,B(x,y)\,z+\overline\omega\,B(y,z)\,x,\label{reduced}\\
     && [x,y,z]^\ast=B(z,x)\,y+\overline\omega\,B(x,y)\,z+\omega\,B(y,z)\,x. \label{reduced conjugate}
\label{reduced form}
\end{eqnarray}
It is easy to verify that the reduced ternary commutator (\ref{reduced}) and its conjugate (\ref{reduced conjugate}) have the same transformation properties under cyclic permutations of arguments (\ref{cyclic}), (\ref{noncyclic}) as the full-length commutator (\ref{ternary commutator}) and its conjugate.
It is interesting to note that in this particular case the reduced ternary commutator (\ref{reduced}) satisfies a reduced version of the basic identity, which contains only ten terms
\begin{equation}
\circlearrowleft \big[[x,y,z],u,v\big]+\circlearrowleft\big[[x,u,y],v,z\big]=0.
\label{reduced identity}
\end{equation}
The basic identity contains two copies of the dihedral group $D_{10}$. The dihedral group $D_{10}$ contains a subgroup of cyclic permutations $\mathbb Z_5$. Thus, the reduced identity (\ref{reduced identity}) is obtained by reducing each copy of the dihedral group $D_{10}$ to its cyclic subgroup $\mathbb Z_5$.

Let us consider a special case of a ternary $\omega$-Lie algebra constructed using ternary multiplication (\ref{ternary product of vectors}). Let us consider $n$-dimensional vectors of $\mathbb C^n$ as row matrices.  Then we can put $B(x,y)=x\,y^T$, where $y^T$ is the column matrix.
Thus we have the ternary Lie algebra, where the vector space of the ternary $\omega$-Lie algebra is the $n$-dimensional complex vector space $\mathbb C^n$ and the ternary commutator is defined by the formula
\begin{equation}
[x,y,z]=z\,x^T\,y+\omega\;x\,y^T\,z+\overline\omega\;y\,z^T\,x.
\label{ternary commutator for vectors}
\end{equation}
In this particular case we can easily compute the structure constants of the ternary $\omega$-Lie algebra. Indeed let $e_1,e_2,\ldots,e_n$ be the canonical basis for $\mathbb C^n$, that is, the $i$th coordinate of a vector $e_i$ is 1, all other coordinates are equal to zero. Then the structure constants of this ternary $\omega$-Lie algebra are
\begin{equation}
C^m_{ijk}=\delta_{ki}\,\delta^m_j+\omega\;\delta_{ij}\;\delta^m_k+\overline\omega\;\delta_{jk}\,\delta^m_i.
\end{equation}
If we calculate the structure constants of the ternary $\omega$-Lie algebra (\ref{ternary commutator for vectors}) for the simplest case $n=2$ then we get
\begin{equation}
[e_1,e_2,e_1]=e_2,\;\;[e_2,e_1,e_2]=e_1.
\end{equation}
Thus we have obtained a realization of the ternary $\omega$-Lie algebra ${\cal L}_2$ introduced at the end of the previous section using vectors of the complex plane and ternary multiplication (\ref{ternary commutator for vectors}).

Let $M_n(\mathbb C)$ be a vector space of complex $n$th order square matrices. Then the ternary product (\ref{ternary product of vectors}) can be applied to $M_n(\mathbb C)$ if we take $B(\Phi,\Psi)=\mbox{Tr}\,(\Phi\,\Psi)$, where $\Phi,\Psi\in M_n(\mathbb C)$. Then the ternary commutator (\ref{reduced}) takes on the form
\begin{equation}
[\Phi,\Psi,\Omega]=\mbox{Tr}\,(\Omega\,\Phi)\;\Psi+\omega\;\mbox{Tr}\,(\Phi\,\Psi)\;\Omega+\overline\omega\;\mbox{Tr}\,(\Psi\,\Omega)\;\Phi.
\label{ternary for square matrices}
\end{equation}
Hence the ternary commutator (\ref{ternary for square matrices}) induces a structure of ternary $\omega$-Lie algebra on a complex vector space $M_n(\mathbb C)$. It is interesting to note that the ternary commutator, which is also constructed by means of the trace and cyclic permutations of arguments
\begin{equation}
\llbracket \Phi,\Psi,\Omega\rrbracket=\mbox{Tr}\,(\Phi)\;{[}\Psi,\Omega]+
                         \mbox{Tr}\,(\Psi)\;{[}\Omega,\Phi]+
                         \mbox{Tr}\,(\Omega)\;{[}\Phi,\Psi],
\label{ternary commutator 3-Lie algebra}
\end{equation}
where square brackets on the right-hand side of this formula stand for commutator of two matrices, that is, $[\Phi,\Psi]=\Phi\,\Psi-\Psi\,\Phi$, induces a structure of 3-Lie algebra on a vector space $M_n(\mathbb C)$. The ternary commutator (\ref{ternary commutator 3-Lie algebra}) was introduced in \cite{Awata-Minich} in order to construct a quantization for generalized Hamiltonian mechanics proposed by Nambu. It should be mentioned that the ternary commutators (\ref{ternary for square matrices}) and (\ref{ternary commutator 3-Lie algebra}) have different properties with respect to permutations of arguments, that is, our ternary commutator (\ref{ternary for square matrices}) has $\omega$-symmetry, while the ternary commutator  (\ref{ternary commutator 3-Lie algebra}) is totally skew-symmetric.

The example of a ternary Lie algebra with ternary commutator (\ref{ternary commutator for vectors}) is a special case of a more general construction. In other words, we can extend the ternary commutator (\ref{ternary commutator for vectors}) to rectangular matrices of arbitrary dimensions. Let $M_{m,n}(\mathbb C)$ be a vector space of complex $m\times n$-matrices. One can define the ternary product of three $m\times n$-matrices $A,B,C\in M_{m,n}(\mathbb C)$ as follows
$$
A\cdot B\cdot C=A\,B^T\,C,
$$
where on the right side of this formula we mean the usual matrix multiplication and $B^T$ stands for transposed matrix. It is easy to verify that this ternary product of $m\times n$-matrices has the associativity of the second kind. Hence we can endow the complex vector space $M_{m,n}(\mathbb C)$ with the following ternary commutator
\begin{equation}
[A,B,C]=A\,B^T\,C+\omega\;B\,C^T\,A+\overline\omega\;C\,A^T\,B+
           C\,B^T\,A+\overline\omega\;B\,A^T\,C+\omega\;A\,C^T\,B,
\end{equation}
and the complex vector space $M_{m,n}(\mathbb C)$ of rectangular $m\times n$-matrices becomes a ternary $\omega$-Lie algebra.

A large class of ternary $\omega$-Lie algebras can be constructed by means of three-dimensional matrices. According to the theorem proposed in \cite{Abramov-Kerner-Liivapuu-Shitov}, there is no ternary product of three-dimensional matrices that satisfies associativity of the first kind. However, the situation with associativity of the second kind of ternary multiplication of three-dimensional matrices is much better. In the same paper \cite{Abramov-Kerner-Liivapuu-Shitov}, the authors found four different ternary products of three-dimensional matrices with associativity of the second kind.
\begin{theorem}
Let $A,B,C$ be $N$th order complex three-dimensional matrices. Then there are only four different triple products of $N$th order complex  three-dimensional matrices which obey the associativity of the second kind. These are
\begin{enumerate}
\item[1)]
$(A\odot B\odot C)_{ijk} = A_{ilm}B_{nlm}C_{njk},\quad
      A\odot B\odot C\rightarrow
\xy <1cm,0cm>:
(1,0)*+{A} , (2,0)*+{B} , (3,0)*+{C} ,
(1.25,-0.2)*+{\bullet} , (1.45,-0.2)*+{\circ} , (1.65,-0.2)*+{\circ} ,
(2.25,-0.2)*+{\circ} , (2.45,-0.2)*+{\circ} , (2.65,-0.2)*+{\circ} ,
(3.25,-0.2)*+{\circ} , (3.45,-0.2)*+{\bullet} , (3.65,-0.2)*+{\bullet} ,
(1.45,-0.26);(2.45,-0.26)**\crv{(1.5, -0.5)&(1.95, -0.7)&(2.4, -0.5)} ,
(1.65,-0.26);(2.65,-0.26)**\crv{(1.7, -0.4)&(2.15, -0.6)&(2.6, -0.4)} ,
(2.25,-0.26);(3.25,-0.26)**\crv{(2.3, -0.5)&(2.75, -0.7)&(3.2, -0.5)}
\endxy$
\item[2)]
$(A\odot B\odot C)_{ijk} = A_{ilm}B_{nml}C_{njk},\quad
     A\odot B\odot C \rightarrow
\xy <1cm,0cm>:
(1,0)*+{A} , (2,0)*+{B} , (3,0)*+{C} ,
(1.25,-0.2)*+{\bullet} , (1.45,-0.2)*+{\circ} , (1.65,-0.2)*+{\circ} ,
(2.25,-0.2)*+{\circ} , (2.45,-0.2)*+{\circ} , (2.65,-0.2)*+{\circ} ,
(3.25,-0.2)*+{\circ} , (3.45,-0.2)*+{\bullet} , (3.65,-0.2)*+{\bullet} ,
(1.45,-0.26);(2.65,-0.26)**\crv{(1.5, -0.5)&(2.05, -0.7)&(2.6, -0.5)} ,
(1.65,-0.26);(2.45,-0.26)**\crv{(1.7, -0.4)&(2.05, -0.6)&(2.4, -0.4)} ,
(2.25,-0.26);(3.25,-0.26)**\crv{(2.3, -0.5)&(2.75, -0.7)&(3.2, -0.5)}
\endxy$
\item[3)]
$(A\odot B\odot C)_{ijk} = A_{ijl}B_{nml}C_{mnk},\quad
    A\odot B\odot C \rightarrow
\xy <1cm,0cm>:
(1,0)*+{A} , (2,0)*+{B} , (3,0)*+{C} ,
(1.25,-0.2)*+{\bullet} , (1.45,-0.2)*+{\bullet} , (1.65,-0.2)*+{\circ} ,
(2.25,-0.2)*+{\circ} , (2.45,-0.2)*+{\circ} , (2.65,-0.2)*+{\circ} ,
(3.25,-0.2)*+{\circ} , (3.45,-0.2)*+{\circ} , (3.65,-0.2)*+{\bullet} ,
(1.65,-0.26);(2.65,-0.26)**\crv{(1.7, -0.5)&(2.15, -0.7)&(2.6, -0.5)} ,
(2.25,-0.26);(3.45,-0.26)**\crv{(2.3, -0.5)&(2.85, -0.7)&(3.4, -0.5)} ,
(2.45,-0.26);(3.25,-0.26)**\crv{(2.45, -0.4)&(2.85, -0.6)&(3.2, -0.4)} ,
\endxy$
\item[4)]
$(A\odot B\odot C)_{ijk} = A_{ijl}B_{mnl}C_{mnk},\quad
      A\odot B\odot C\rightarrow
\xy <1cm,0cm>:
(1,0)*+{A} , (2,0)*+{B} , (3,0)*+{C} ,
(1.25,-0.2)*+{\bullet} , (1.45,-0.2)*+{\bullet} , (1.65,-0.2)*+{\circ} ,
(2.25,-0.2)*+{\circ} , (2.45,-0.2)*+{\circ} , (2.65,-0.2)*+{\circ} ,
(3.25,-0.2)*+{\circ} , (3.45,-0.2)*+{\circ} , (3.65,-0.2)*+{\bullet} ,
(1.65,-0.26);(2.65,-0.26)**\crv{(1.7, -0.5)&(2.15, -0.7)&(2.6, -0.5)} ,
(2.25,-0.26);(3.25,-0.26)**\crv{(2.3, -0.5)&(2.75, -0.7)&(3.2, -0.5)} ,
(2.45,-0.26);(3.45,-0.26)**\crv{(2.45, -0.4)&(2.95, -0.6)&(3.4, -0.4)} ,
\endxy$
\label{theorem_ternary_multiplications}
\end{enumerate}
In the diagrammatic representation of ternary multiplication, one should take a sum over a pair of indices depicted by empty circles connected by arcs and black filled circles represent free indices.
\end{theorem}
We will consider the simplest example of a ternary $\omega$-Lie algebra constructed using three-dimensional matrices of the second order. As a ternary product of three-dimensional matrices, we will use ternary multiplication 3 (Theorem \ref{theorem_ternary_multiplications}), although it is worth noting that we could equally well use ternary multiplication 4. Let $A=(a_{ijk})$ be a three-dimensional matrix of second order, that is, $i,j,k=1,2$. We will call a three-dimensional matrix $A$ {\em traceless} if the trace of this matrix with respect to any pair of subscripts is zero. Hence for any $k=1,2$ we have
\begin{equation}
a_{iik}=a_{iki}=a_{kii}=0,
\label{traceless}
\end{equation}
where $a_{iik}=a_{11k}+a_{22k}, a_{iki}=a_{1k1}+a_{2k2},a_{kii}=a_{k11}+a_{k22}.$ The ternary $\omega$-Lie algebra of three-dimensional matrices of the second order is an 8-dimensional algebra. Traceless matrices form a two-dimensional subspace in this algebra, and it is easy to show that this two-dimensional subspace is closed under the ternary commutator (\ref{ternary commutator}), that is, traceless three-dimensional matrices of the second order form a subalgebra of the ternary $\omega$-Lie algebra of three-dimensional matrices of the second order. From conditions (\ref{traceless}) it follows that in the case of a three-dimensional matrix of the second order $A$ we have two independent parameters $a_{111},a_{222}$, and all other entries of the matrix are expressed through them, that is,
$$
a_{221}=a_{212}=a_{122}=-a_{111},\;\;a_{112}=a_{121}=a_{211}=-a_{222}.
$$
We arrange the entries of a three-dimensional matrix of 2nd order $A$ in space, that is, in the vertices of the cube, as follows
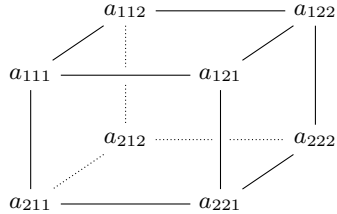
\begin{center}

\begin{tikzpicture}
  \matrix (m1) [matrix of math nodes, row sep=1.5em, column sep=1.5em]{
    & a_{112} & & a_{122} \\
    a_{111} & & a_{121} & \\
    & a_{212} & & a_{222} \\
    a_{211} & & a_{221} & \\};

  \path[-]
    (m1-1-2) edge (m1-1-4)
            edge (m1-2-1)
            edge [densely dotted] (m1-3-2)
    (m1-1-4) edge (m1-3-4)
            edge (m1-2-3)
    (m1-2-1) edge [-,line width=6pt,draw=white] (m1-2-3)
            edge (m1-2-3)
            edge (m1-4-1)
    (m1-3-2) edge [densely dotted] (m1-3-4)
            edge [densely dotted] (m1-4-1)
    (m1-4-1) edge (m1-4-3)
    (m1-3-4) edge (m1-4-3)
    (m1-2-3) edge [-,line width=6pt,draw=white] (m1-4-3)
            edge (m1-4-3);
 \end{tikzpicture}
  \end{center}
Thus, as generators of the ternary $\omega$-Lie algebra of three-dimensional traceless matrices of the second order, we can take two three-dimensional matrices $F_1=-\frac{i}{2\sqrt{2}}\,E_1,F_2= -\frac{i}{2\sqrt{2}}\,E_2$, where

\begin{tikzpicture}
  \matrix (m1) [matrix of math nodes, row sep=1.5em, column sep=1.5em]{
    & 0 & & -1 \\
    1 & & 0 & \\
    & -1 & & 0 \\
    0 & & -1 & \\};
   \node[left=0.1em of m1] {\(E_1=\)};
  \path[-]
    (m1-1-2) edge (m1-1-4)
            edge (m1-2-1)
            edge [densely dotted] (m1-3-2)
    (m1-1-4) edge (m1-3-4)
            edge (m1-2-3)
    (m1-2-1) edge [-,line width=6pt,draw=white] (m1-2-3)
            edge (m1-2-3)
            edge (m1-4-1)
    (m1-3-2) edge [densely dotted] (m1-3-4)
            edge [densely dotted] (m1-4-1)
    (m1-4-1) edge (m1-4-3)
    (m1-3-4) edge (m1-4-3)
    (m1-2-3) edge [-,line width=6pt,draw=white] (m1-4-3)
            edge (m1-4-3);

  \matrix (m2) [matrix of math nodes, row sep=1.5em, column sep=1.5em, right=5em of m1]{
    & -1 & & 0 \\
    0 & & -1 & \\
    & 0 & & 1 \\
    -1 & & 0 & \\};
  \node[left=0.1em of m2] {\(E_2=\)};
  \path[-]
    (m2-1-2) edge (m2-1-4)
            edge (m2-2-1)
            edge [densely dotted] (m2-3-2)
    (m2-1-4) edge (m2-3-4)
            edge (m2-2-3)
    (m2-2-1) edge [-,line width=6pt,draw=white] (m2-2-3)
            edge (m2-2-3)
            edge (m2-4-1)
    (m2-3-2) edge [densely dotted] (m2-3-4)
            edge [densely dotted] (m2-4-1)
    (m2-4-1) edge (m2-4-3)
    (m2-3-4) edge (m2-4-3)
    (m2-2-3) edge [-,line width=6pt,draw=white] (m2-4-3)
            edge (m2-4-3);
\end{tikzpicture}

\noindent
Calculating the ternary commutator
\begin{eqnarray}
[A,B,C]_{ijk}\!\!\! &=&\!\!\! A_{ijl}B_{nml}C_{mnk}+\omega\;B_{ijl}C_{nml}A_{mnk}+
      \overline\omega\;C_{ijl}A_{nml}B_{mnk}\nonumber\\
      &&\quad +C_{ijl}B_{nml}A_{mnk}+
      \overline\omega\;B_{ijl}A_{nml}C_{mnk}+\omega\;A_{ijl}C_{nml}B_{mnk},\nonumber
\end{eqnarray}
we find the commutation relations of the ternary $\omega$-Lie algebra of three-dimensional traceless matrices of the second order
$$
[F_1,F_2,F_1]=F_2,\;\;[F_2,F_1,F_2]=F_1.
$$
Thus we constructed one more realization of the 2-dimensional ternary $\omega$-Lie algebra ${\cal L}_2$ by means of traceless three-dimensional matrices of second order.

\section{Discussion}
In this paper we propose an extension of the concept of Lie algebra to algebras with ternary multiplication laws. Our approach is based on the concept of a ternary commutator, which we construct by analogy with a binary commutator, that is, we form six triple products using all permutations of three arguments of ternary commutator and then form their linear combination using the 3rd order roots of unity as coefficients. Due to the properties of 3rd roots of unity, the ternary commutator we propose vanishes identically when all three of its arguments are equal. However, in the case when two of its three arguments are equal, it is generally not equal to zero.
Here we note the analogy with the ternary generalization of the Pauli exclusion principle proposed by Kerner \cite{Kerner(2017)}. According to this generalization, a wave function of a quantum system consisting of three particles vanishes identically when all three particles have identical quantum characteristics. However, in the case where two particles have equal quantum characteristics, a wave function does not necessarily vanish. As an example, we can point to the properties of the quark model. The theory of groups, Lie algebras and their representations is successfully used in the theory of elementary particles. We think that the concept of a ternary Lie algebra at cube root of unity, introduced in this paper, adequately reflects the basic properties of the quark model.

\vspace{6pt}


\begin{thebibliography}{999}
\bibitem{Abramov-Kerner-LeRoy}
Abramov, V., Kerner, R. and Le Roy, B. {Hypersymmetry: A $\mathbb Z_3$-generalization of supersymmetry}, {\em J. Math. Phys.} {\bf 1997}, {\em 38} (3), 1650--1669.
\bibitem{Abramov-Kerner-Liivapuu-Shitov}
Abramov, V., Kerner, R., Shitov, S. and Liivapuu, O. Algebras with ternary law of composition and their realization by cubic matrices, {\em Jornal of Generalized Lie Theory and Applications} {\bf 2009}, {\em 3} (2), 77--94.
\bibitem{Abramov-Liivapuu}
Abramov, V. and Liivapuu, O. SO(3)-Irreducible Geometry in Complex Dimension Five and Ternary Generalization of Pauli Exclusion Principle, {\em Universe} {\bf 2024}, {\em 10} (1), 1.
\bibitem{Awata-Minich}
Awata, H., Li, M., Minic, D., and Yaneya, T.
{On the quantization of Nambu brackets}, {\em JHEP02} {\bf 2001}, {\em 013}.
\bibitem{Bagger-Lambert(2007)}
Bagger, J. and Lambert, N. {Modeling multiple M2's}. {\em Phys. Rev} {\bf 2007}, {\em D75}, 045020, \href{https://arxiv.org/abs/hep-th/0611108}{arXiv:hep-th/0611108}.
\bibitem{Bagger-Lambert(2008)}
Bagger, J. and Lambert, N. {Gauge symmetry and supersymmetry of multiple M2-branes}. {\em Phys. Rev} {\bf 2008}, {\em D77}, 065008, \href{https://arxiv.org/abs/0711.0955}{arXiv:0711.0955 [hep-th]}.
\bibitem{Filippov(1985)}
Filippov, V. T. {$n$-Lie algebras}. {\em Siberian Math. J.} {\bf 1985}, {\em 26}, 879--891.
\bibitem{Gelfand-Minlos-Shapiro}
Gelfand, I.~M., Minlos, R.~A., Shapiro, Z.~Ya. {Representations of the Rotation and Lorentz Groups and Their Applications}, Dover Publications, Ins. Mineola, New York, 2018.
\bibitem{Kerner(2017)}
Kerner, R. {Ternary Generalization of Pauli's Principle and the $Z_6$-Graded Algebras}, {\em Phys. At. Nucl.} {\bf 2017}, {\em 80}, 522--534.
\bibitem{Kerner(2019)}
Kerner, R. {The Quantum Nature of Lorentz Invariance}, {\em Universe} {\bf 2018}, {\em 5} (1).
\bibitem{Nambu(1973)}
Nambu, Y. {Generalized Hamiltonian mechanics}. {\em Phys. Rev. D} {\bf 1973}, {\em 7}, 2405--2412.
\bibitem{Wagner(1953)}
Wagner, V.~V. The theory of generalized heaps and generalized groups, {\em Matematicheskii Sbornik} {\bf 1953}, {\em 74} (3), 545 -- 632.
\bibitem{Zapata-Arsiwalla-Beynon(2024)}
Zapata-Carratal\'{a}, C., Arsiwalla, X.~D., Beynon, T. Heaps of Fish: arrays, generalized associativity and heapoids, \href{https://arxiv.org/abs/2205.05456}{arXiv:2205.05456[math.RA]} (to appear in {\em Theoretical Computer Science}).
\end{thebibliography}
\end{document}